\documentclass[11pt]{amsart}
\usepackage{graphicx}
\usepackage{amssymb}
\usepackage{color}

\newcommand{\R}{\mathbb{R}}

\def\C{{\mathbb C}} 
\def\R{{{\mathbb R}}}

\newcommand{\norm}[1]{\left\Vert#1\right\Vert}
\newcommand{\abs}[1]{\left\vert#1\right\vert}
\newcommand{\ip}[1]{\langle#1\rangle}
\def\be{\begin{equation}}
\def\ee{\end{equation}}
\def\eps{{\epsilon}}

\newtheorem{theorem}{Theorem}

\newtheorem{remark}[theorem]{Remark}

\newtheorem{lemma}[theorem]{Lemma}

\numberwithin{equation}{section}
\numberwithin{theorem}{section}

\title[Blow-up rate of solutions
 of the  3D Zakharov system] {  Lower bound  for 
the rate of blow-up of singular solutions
 of the Zakharov system in $\R^3$}
\author{J.~Colliander}
\address{Department of Mathematics, University of Toronto}

\thanks{J.C. is partially supported by NSERC through grant number RGPIN 250233-12.}
\author{M. ~Czubak}
\address{Department of Mathematical Sciences, Binghamton University (SUNY)}
\author{C.~Sulem}
\address{Department of Mathematics, University of Toronto}

 \thanks{C.S. is  partially  supported by NSERC through grant number RGPIN  46179-11.}

\subjclass{35Q55}
\keywords{Zakharov system, local wellposedness, singular solutions, 
blow-up rate}

\begin{document}
\date{\today}

\begin{abstract}
We consider the scalar Zakharov system in $\R^3$ for initial conditions
$(\psi(0), n(0), n_t(0)) \in  H^{\ell+1/2} \times H^\ell \times H^{\ell-1}  $,
 $0\leq\ell \leq 1$.
Assuming that the solution blows up in a finite time $t^* < \infty$, 
we establish a lower bound
for the  rate of blow-up  of the corresponding Sobolev norms in the form
$$ \|\psi(t)\|_{H^{\ell+1/2}} +\|n(t)\|_{H^{\ell}} + \|n_t(t)\|_{H^{\ell-1}}
> C(t^*-t)^{-\theta_\ell} $$
with $\theta_\ell = \frac{1}{4}(1+ 2  \ell)^-$. 
The analysis is a reappraisal of the local wellposedness theory of 
Ginibre, Tsutsumi and Velo (1997) combined with an argument 
developed by  Cazenave and Weissler (1990) in the context of 
nonlinear Schr\"odinger equations.

\end{abstract}

\maketitle

\section{Introduction}

The Zakharov  system describes the phenomenon of propagation of Langmuir
waves in a non-magnetized plasma. It was derived by Zakharov \cite{Z72}
in the form of  a coupled system governing the electric field  complex amplitude $\psi(x,t)$
and the density fluctuations of ions $n(x,t)$. 
Here we consider the scalar   Zakharov  system in the form:
\begin{eqnarray}
&& i\partial_t \psi + \Delta \psi = n  \psi ,  \label{zakh1}\\
&&\partial_{tt}n - \Delta n = \Delta  |\psi|^2, \label{zakh2}
\end{eqnarray}
where $\psi : (x,t) \in \R^d \times \R^+ \to \C$, $n: (x,t)\in \R^d\times \R^+ \to \R$, 
with  given initial conditions $\psi(0)=\psi_0$, $n(0)=n_0$ and $n_t(0)=n_1$.

There has been a large literature devoted to the local and global wellposedness
of the initial value problem in the context of smooth solutions
(\cite{SS79}, \cite{AA84}, \cite{SW86}, \cite{OT91},
\cite{GM94}). In  recent years, an effort has been made to lower the
regularity assumptions (\cite{BC96}, \cite{GTV97}, \cite{P01}, \cite{CHT08},
\cite{BHHT09}, \cite{BH11}, \cite{T00}) and  to investigate the possible occurrence of 
 local ill-posedness \cite{H07}.

For general initial conditions in the energy space with negative Hamiltonian, solutions to the Zakharov system
in two and three dimensions will blow up in a finite or infinite  time \cite{M96b}. Heuristic arguments
and numerical simulations show that solutions do  blow up
in a finite time  both in two and three dimensions (see \cite{SS99} for a review).

In two dimensions, there exist exact self-similar blow-up solutions 
\cite{ZS81}
\begin{eqnarray}
&&\psi( x,t)= \frac{1}{a (t^*-t)} P\Big(\frac{|x|}{a(t^* -t)}
\Big)
e^{i\left (\theta + \frac{1}{a^2(t^*-t)}-\frac{|x|^2}{4(t^*-t)}\right )}, \label{selfsimzak1}\\
&&n(x,t) = \frac{1}{a^2(t^*-t)^2} N\Big(\frac{|x|}{a (t^*-t)}\Big),
\label{selfsimzak2}
\end{eqnarray}
where $(P,N$) satisfy the system of ODEs in the radial variable denoted  $\eta$
\begin{eqnarray}
&&\Delta P-P-NP =0, \label{Ppeq2} \\
&&a^2(\eta^2 N_{\eta \eta}+6\eta N_{\eta}+6N)- \Delta N=\Delta P^2,
\label{Nneq2}
\end{eqnarray}
and $a>0$ is a free parameter.
Rigorous results on these solutions were proved in \cite{GM94}.
Numerical simulations show that for a large class of data, blow-up
solutions   asymptotically  display a self-similar collapse  described by  the above solutions
\eqref{selfsimzak1}-\eqref{selfsimzak2}.
In addition, Merle \cite{M96a} established  a lower bound for the rate of blow-up 
of singular  solutions of the Zakharov system  in the energy space
in the form
\begin{equation}
\|u(t)\|_{H^{1}} \ge C(t^*-t)^{-1} \, , \, \, \|n(t)\|_{L^2} 
\ge C(t^*-t)^{-1}.
\end{equation}
This rate is optimal, in the sense that the exact solutions 
\eqref{selfsimzak1}-\eqref{selfsimzak2}
solutions blow up exactly in this fashion.  It is an open question whether  there exist other solutions
that blow up at a faster rate.

 Merle  uses a time-dependent rescaling based on
the scale invariance of the wave equation.  The scaling factor is associated to the energy norm
 and the energy conservation is interpreted in
terms of the new variables. The optimal constant for Sobolev inequality  expressed in terms
of the ground state of the 2d cubic NLS equation
is used to obtain a lower bound for the scaling factor, which in turn is related
to the energy norm of the solution. A completely new element in \cite{M96a} was a compactness argument leading to a limiting  object as $t$ approaches $t^*$. This method, now referred to  as `the Liouville approach' opened doors to
break-through developments in the field.

The situation for the Zakharov system in three dimensions is more complex
 and several questions remain open. There are no known explicit
 blow-up solutions.
Self-similar solutions exist only asymptotically close to collapse. 
They have the universal form
\cite{BZS75}, \cite{ZS81}
\begin{eqnarray}
\psi(x,t) &=& \frac{1}{(t^*-t)}P\Big(\frac{|x|}{(t^*-t)^{2/3}}\Big)
e^{i(t^*-t)^{-1/3}}, \label{r.ss}\\ 
n(x,t)  &=& \frac{1}{(t^*-t)^{4/3}}
N\Big(\frac{|x|}{(t^*-t)^{2/3}}\Big),\label{n.ss} 
\end{eqnarray}
where $P(\eta)$ and $N(\eta)$ are radially symmetric scalar functions satisfying 
the coupled system of ODEs
\begin{eqnarray}
\Delta P-\frac{1}{3}P-NP&=&0,  \label{Req} \\
\frac{2}{9}(2\eta^2 N_{\eta\eta} + 13 \eta N_{\eta}
+ 14 N)&=&\Delta {P}^2  \label{Neq}.
\end{eqnarray}

Note that there is no free parameter in this system.  In addition,  unlike the 2d case,
there is no rigorous proof of existence of solutions to the ODE system
 \eqref{Req}-\eqref{Neq}. The profile system \eqref{Req}-\eqref{Neq}
 was studied numerically in \cite{ZS81} where two pairs of localized
solutions were computed; one of them displaying a monotone profile for $P$ and $N$. 
Like for the 2d problem, numerical simulations  of the three-dimensional 
Zakharov system \cite{LPSSW92}
 show that, for a large class of initial conditions,  solutions 
 blow up in a finite time and display a  self-similar
collapse described by \eqref{r.ss}-\eqref{n.ss}. This can be seen 
as the dynamic stability of these  asymptotic solutions.  
Asymptotically close to the collapse, the regime is 
strongly supersonic
with the pressure term
$\Delta n$ negligible compared to
the ion-inertia term in  \eqref{zakh2}.

In the present note, we consider the question of  the rate of blow-up of solutions
for the Zakharov system in three dimensions, and
we establish a lower bound for it in 
appropriate Sobolev norms.
Our method differs from that developed by Merle \cite{M96a} for the 2d
problem.  It is in a sense simpler but less precise. The two main ingredients are a local well-posedness result  and
a contradiction argument adapted from Cazenave and Weissler \cite{CW90}.

The notion of criticality plays a central role in the study of the Nonlinear
Schr\"odinger equation (NLS). For the Zakharov system however, 
criticality is less straightforward because the NLS and the wave equation
scale differently. In \cite{GTV97}, Ginibre, Tsutsumi and Velo 
proposed a definition of criticality  for the Zakharov system 
for initial condition $(\psi(0), n(0), n_t(0))$ in 
$H^k \times H^\ell \times H^{\ell-1}$ with the critical values being
$k=d/2-3/2$ , $\ell=d/2-2$, and $d$ is the spatial dimension.
Note that $k-\ell=1/2$ and not one as the energy space would suggest.

In three dimensions,  the critical space in the above sense is 
$L^2\times H^{-1/2} \times  H^{-1}$ which is, up to $\epsilon>0$ the 
space in which Bejenaru and Herr \cite{BH11} recently proved local
well-posedness. 
Also, the asymptotic solution \eqref{r.ss}-\eqref{n.ss} has the
property that the $H^k$ norm of $\psi$ and the $H^\ell $ norm of $n$ 
blow up at the same rate when $k -\ell =1/2$.

Our analysis  relies 
on the local well-posedness results of Ginibre et al
\cite{GTV97} in  $H_\ell =  H^{\ell+1/2}  \times H^{\ell} 
\times H^{\ell-1}$,
$\ell\ge 0$, thus concerns solutions that are slightly 
more regular than solutions in the critical space.

\begin{theorem}\label{bup}

 Consider 
 initial conditions $(\psi(0), n(0), n_t(0))$
  in  
$H_\ell$,
 $0\leq \ell \leq 1$. Assume that the solution $(u,n)$ blows up in a finite time
$t^*$ in $H_l$.  Then, we have the lower bound for the rate of blow-up in the 
corresponding Sobolev norms 
\begin{equation}
\|\psi(t)\|_{H^{\ell+1/2}} +\|n(t)\|_{H^{\ell}} + \|n_t(t)\|_{H^{\ell-1}}
> C(t^*-t)^{-\theta_\ell}
\end{equation}
with $\theta_\ell = \frac{1}{4}(1+ 2 \ell)^-$.
\end{theorem}

\begin{remark}
We do not know whether this bound is optimal.
 In particular, we observe that
the homogeneous $\dot H^{\ell+1/2}$ norm of $\psi$ and the homogeneous $\dot H^{\ell}$ norm of $n$ in the expression
of the 
asymptotic solution \eqref{r.ss}-\eqref{n.ss} both blow up at  a faster rate, namely $\frac{1}{3}( 1+ 2 \ell)$.   However, as \eqref{r.ss}-\eqref{n.ss} is not a solution of \eqref{zakh1}-\eqref{zakh2}, but only an asymptotic solution,
its rate of blow-up might not be a real threshold. Moreover, for the cubic NLS in 3D, the present method gives a rate of blow-up of the $H^{1}$ norm to be $\frac 14$ \cite{CW90}, which has been observed in numerical
simulations  \cite[Chapter 7]{SS99}.
\end{remark}
\begin{remark}
We choose to consider the $H^{\ell+1/2}\times H^{\ell}\times H^{\ell-1}$ 
functional framework instead of the more general norm
$H^{k}\times H^{\ell}\times H^{\ell-1}$, because for $k=l+\frac 12$, $\norm{\psi(t)}_{\dot H^k}$ scales the same as $\norm{n(t)}_{\dot H^\ell}$ when $\psi$ and $n$ are given by \eqref{r.ss}-\eqref{n.ss}.
\end{remark}
Here is a brief description of the content of the paper.
In Section 2, we recall  important linear estimates. Section 3 is
devoted to nonlinear estimates. In particular, we carefully keep track of the
power of  time  involved in the estimates as it is central for the analysis of the lower bound for
the blow-up rate. In Section 4, we adapt an argument  for semilinear 
heat equations due to  Weissler \cite{W81} and later extended to 
nonlinear Schr\"odinger equations by  Cazenave and Weissler \cite{CW90} to obtain a
lower bound of blow-up for Sobolev norms of the solution.

\section{Preliminary estimates}

Consider the Zakharov system \eqref{zakh1}-\eqref{zakh2}
with initial conditions 
\begin{equation} \label{ic}
(\psi,n,n_t)\vert_{t=0} ~ = ~(\psi_0,n_0, n_1).
\end{equation}
The wave equation \eqref{zakh2} can be transformed  into a reduced wave equation \cite{BHHT09}, \cite{GTV97} for 
$$w =n+ i\ip{\nabla}^{-1} \partial_t n,$$
 where $\ip{\nabla} =(1-\Delta)^{1/2}$. 
 
 The new system then takes the form 
 \begin{align}
i\partial_t \psi + \Delta \psi &=( \mathcal Re~w)  \psi ,  \label{zakh1a}\\
(i\partial_t - \ip{\nabla} ) w &= -\ip{\nabla}^{-1} \Delta|\psi|^2-\ip{\nabla}^{-1}\mathcal Re~w \label{zakh2a},
\end{align}
and $(\psi, w)$ solve \eqref{zakh1a}-\eqref{zakh2a} with data $(\psi_{0}, w_{0})=(\psi_{0}, n_{0}+i\ip{\nabla}^{-1}n_{1})$ if and only if $(\psi, \mathcal Re~w)$ solve \eqref{zakh1}-\eqref{zakh2} with data $(\psi_{0}, n_{0}, n_{1})$.

We will  use space-time norms in the context of
solutions defined on a finite time interval $(-T,T)$, and we introduce
an even time cut-off function $\varphi \in C^\infty$  satisfying
$\varphi(t) =1$ for $|t|\leq 1$,   $\varphi(t) =0$ for $|t|\ge 2$,
 $0\le \varphi(t)\le 1$.  We  denote  $\varphi_T(t) =\varphi(t/T)$,
($T\le1$).
The Duhamel representation of the solution takes the form
%
\begin{align}
\psi(t) &=\varphi_1(t) U(t) \psi_0-i\varphi_T(t)\int_0^t U(t-s) f_{1}(s) ds,\label{d4psi}\\
 w(t)&= \varphi_1(t) W(t) w_{0}+ i\varphi_T(t)\int_0^t W(t-s)\left(f(s)+\varphi_{2T}\tfrac{\mathcal Re~w}{\ip{\nabla}}\right) ds,\label{d4n}
\end{align}
where
\begin{align*}
 &U(t) = e^{it \Delta},\quad W(t) =e^{- it \sqrt{1-\Delta}},\\
&f_{1}=\varphi_{2T}^2 ~(\mathcal Re~w) \psi, \quad
f=\ip{\nabla}^{-1}\varphi_{2T}^2\Delta|\psi|^2.
\end{align*}

Building on the foundation established in \cite{B93} and following \cite{GTV97}, we seek a solution
 $(\psi, w)\in X_S^{l+\frac 12,b}\times X_{W}^{l,b}$, which are 
the space-time weighted Bourgain spaces with norms respectively given by
\begin{align*}
  \norm{\psi}_{X_S^{l+\frac 12,b}}&=\norm{\ip{\xi}^{l+\frac 12}\ip{\tau+\abs{\xi}^2}^{b}\hat\psi(\tau,\xi)}_{L^2_{\tau,\xi}},\\
\norm{w}_{X_{W}^{l,b}}&=\norm{\ip{\xi}^{l}\ip{\tau+\abs{\xi}}^{b}\hat w(\tau,\xi)}_{L^2_{\tau,\xi}},
\end{align*}
where we  use the notation $\ip{\xi} = (1 +|\xi|^2)^{1/2}$.
Note the difference in the dispersive weights for the above two norms.  We are using $\ip{\tau+\abs{\xi}}$ for the reduced wave equation, which is equivalent to $\ip{\tau+\ip{\xi}}.$  Also, we did not find a benefit of using two different $b$ indices. 

We now recall important linear estimates from \cite{GTV97} (see also \cite{BHHT09}).
\begin{lemma}\label{l:21}Let $s, b\in \R$ and $(X^{s,b}, {\bf\tilde U(t)})=(X^{s,b}_{S},  U(t))$ or $(X^{s,b}, {\bf\tilde U(t)})=(X^{s,b}_{W}, W(t)).$  Then
\begin{equation} 
\|\varphi_1 {\bf \tilde U}(t) u_0 \|_{X^{s,b}} = \|\varphi_1\|_{H^b}\|u_0\|_{H^s}.
\end{equation}
Let $-1/2 < b'\le 0 \le b\le b'+1$, and $T\le 1$. Then 
\begin{equation}\label{l21}
\| \varphi_T \int_0^t {\bf\tilde U}(t-t') f(s) dt'\|_{X^{s,b}}
 \le C T^{1-b+b'}\| f\|_{X^{s,b'}}.
\end{equation}
\end{lemma}
The cut-off function $\varphi_{2T}$ has been introduced inside the nonlinear
term in \eqref{d4psi}-\eqref{d4n}. Its effect is evaluated in the next lemma.

\begin{lemma} For any $s\in \R$, $b\ge0$, $q\ge 2$ and $bq >1$, 
\begin{equation}\label{l25}
\|\varphi_T  f \|_{X^{s,b}} \le C T^{-b+1/q} \|f \|_{X^{s,b}},
\end{equation}
where $X^{s,b}=X^{s,b}_{S}$ or $X^{s,b}=X^{s,b}_{W}$.
\end{lemma}
Note that for a  parameter   
 $b>1/2$, the negative power of $T$ is minimized with $q=2$.
  
We apply Lemma \ref{l:21} to \eqref{d4psi} with
\[
b'=b-1+\eps, \ \ 0<  \eps \ll 1,
\]
and obtain
\be\label{dh1}
\begin{split}
\norm{\psi}_{X_S^{\ell+\frac 12,b}}&\lesssim \norm{\psi_0}_{H^{\ell+\frac 12}}+T^{\eps}\norm{\varphi^2_{2T}(\mathcal Re~w) \psi}_{X_S^{\ell+\frac 12,b-1+\eps}}.
\end{split}
\ee
Similarly for \eqref{d4n}, we first apply Lemma \ref{l:21} to the nonlinear term $f$ with $b'=b-1+\eps, \ \ 0<  \eps \ll 1$ and then to the linear term with $b'=0$.  This results
in
\begin{align}
\norm{w}_{X_{W}^{\ell,b}}&\lesssim \norm{w_{0}}_{H^{\ell}}\nonumber\\
&\quad+T^{\eps}\norm{\varphi^2_{2T}\frac{\Delta}{\ip{\nabla}}\abs{\psi}^2}_{X_{W}^{\ell,b-1+\eps}}+
T^{1-b}\norm{\varphi_{2T}  \ip{\nabla}^{-1} \mathcal Re~w}_{X_{W}^{\ell,0}}
.\label{dh2}
\end{align}
The next section is dedicated to showing we can handle the nonlinearities on the right hand side, and produce \emph{additional} powers of $T$ in the process.


\section{Nonlinear estimates}
We need to estimate the right hand side of \eqref{dh1} and \eqref{dh2}, namely
\begin{align*}
f_{1}=\varphi^2_{2T}( \mathcal Re~w) \psi \in X_S^{k,-c},\quad
f= \ip{\nabla}^{-1}\varphi_{2T}^2\Delta|\psi|^2 \in X_{W}^{\ell, -c},
\end{align*}
where $c=-(b-1+\eps),$ and $k-\ell=\frac 12$. The second term in \eqref{dh2} will be treated in the next section. 
More precisely, we need to establish
\begin{align} 
\norm{\varphi^2_{2T}( \mathcal Re~w)\psi}_{X_S^{k,-c}}&\lesssim 
T^{\theta}\norm{\varphi_{2T} \mathcal Re~w}_{X_{W}^{\ell,b}}\norm{\varphi_{2T}\psi}_{X_S^{k,b}}, \label{m10}\\
\norm{   \ip{\nabla}^{-1}  \left(\varphi_{2T}^2\Delta|\psi|^2\right) }_{X_{W}^{\ell,-c}}
&\lesssim T^{\theta}
\norm{\varphi_{2T} \psi}_{X_S^{k,b}}^{2}.\label{m20}
\end{align}

The following setup is standard.  Let 
$\varphi_{2T}\mathcal Re~w=u$.  We consider the nonlinearities on the Fourier side following the notation of \cite{GTV97}.  
\begin{eqnarray}
&&\hat f_1(\xi_1,\tau_1)=\int_{\R^{3+1}}\widehat{u}(\xi_1-\xi_2,\tau_1-\tau_2)
(\widehat{\varphi_{2T} \psi})(\xi_2,\tau_2) d\xi_2d\tau_2,\\
&&|\hat f(\xi,\tau)|\leq\abs{\xi}\int_{\R^{3+1}}\widehat{\abs{\varphi_{2T}\psi}}(\xi+\xi_2,\tau+\tau_2)
\widehat{ \abs{\varphi_{2T}\psi}}(-\xi_2,-\tau_2) d\xi_2d\tau_2,
\end{eqnarray}
where we changed variables in the second integral, and used the trivial estimate $\abs{\xi}^2\ip{\xi}^{-1}\leq \abs{\xi}$. 

To estimate $f_1$ in $X^{k,-c}_{S}$, we  define
 \begin{eqnarray}
 &&\hat v_2(\xi_{2,},\tau_{2}) = \ip{\xi_2}^k \ip{ \tau_2 +\abs{\xi_2}^2}^{b} (\widehat{\varphi_{2T} \psi}) (\xi_2, \tau_2),\\
 &&\hat v(\xi,\tau) = \ip{\xi}^l \ip{\tau +\abs{\xi}}^b\hat u(\xi,\tau),
 \end{eqnarray}
 and by duality,  take the scalar product with a test function
 in $X^{-k,c}_{S}$ or equivalently with a function whose Fourier transform 
 is $ \ip{\xi_1}^k \ip{\tau_1 +\abs{\xi_1}^2}^{-c}\hat v_1(\xi_1,\tau_1) $ with 
 $v_1$ in $L^2_{xt}$. 
 
 Proceeding similarly for $f$,
 we are led to showing two estimates
\begin{align}
\abs{N_1(v, v_1, v_2)}\lesssim T^\theta\norm{v}_2\norm{v_1}_2\norm{v_2}_2\label{m1},\\
\abs{N_2(v, v_1, v_2)}\lesssim T^\theta\norm{v}_2\norm{v_1}_2\norm{v_2}_2,\label{m2}
\end{align}
where

\begin{align}
N_1(v, v_1, v_2)&=\int\frac{\hat v \hat v_1\hat v_2~\ip{\xi_1}^k  ~d\xi_1 d\xi_2 d\tau_1 d\tau_2}
{\ip{\tau+\abs{\xi}}^b\ip{\tau_1+\abs{\xi_1}^2}^{c}\ip{\tau_2+\abs{\xi_2}^2}^{b}\ip{\xi_2}^k\ip{\xi}^l } ,\\
N_2(v, v_1, v_2)&=\int\frac{\hat v \hat v_1\hat v_2~\abs{\xi}~\ip{\xi}^\ell ~d\xi_1 d\xi_2 d\tau_1 d\tau_2 }
{\ip{\tau+\abs{\xi}}^c\ip{\tau_1+\abs{\xi_1}^2}^{b}\ip{\tau_2+\abs{\xi_2}^2}^{b}\ip{\xi_1}^k\ip{\xi_2}^k }.
\end{align}
The arguments of $\hat v $,  $\hat v_1 $,  $\hat v_2 $  are $(\xi,\tau)$,
$(\xi_1,\tau_1)$, $(\xi_2,\tau_2)$ with  $\xi = \xi_{1}-\xi_2$, $\tau = \tau_{1}- \tau_2$, and
we can also assume $\mathcal F^{-1} (\tfrac{ \hat v}{\ip {\tau+\abs{\xi}}^{\gamma a}}), \mathcal F^{-1} (\tfrac { \hat v_i}{\langle \tau_i+\abs{\xi_i}^2 \rangle^{\gamma a_i}}), i=1, 2, $ have support in $\abs{t}\leq CT$.

Introducing the variables $\sigma_i= \tau_i+\abs{\xi_i}^2$ and $\sigma= \tau+\abs{\xi}$, we can
express
\[
\xi_1^2-\xi^2_2-\abs{\xi}=\sigma_1-\sigma_2-\sigma,
\]
from which one  concludes (using ideas first observed in \cite{B94} and used in \cite{BC96} and \cite[Lemma 3.3]{GTV97} )
\be\label{bound}
\ip{\xi_1}^2\lesssim \ip{\sigma_1}+\ip{\sigma_2}+\ip{\sigma}\quad \mbox{when} \quad \abs{\xi_1}\geq 2\abs{\xi_2}.
\ee
In \cite{GTV97}, the proof of the two estimates in the range of exponents we are interested in, 
was accomplished in two lemmas (Lemma 3.4 and 3.5), which were obtained from a repeated
application of a general estimate (shown in Lemma 3.2).  The analysis of \cite{GTV97} 
did not require an optimal power of $\theta$, but needed it to be just large enough, so the 
final power of $T$, after combining all the estimates, was positive.  Hence it appears that certain simplifications were made, which resulted in cleaner estimates (See for example 
\cite[Remark 3.1]{GTV97}.  Here we seek the optimal power of 
$\theta$, so we reprove estimates \eqref{m1}-\eqref{m2} with a goal 
to optimize the final power of $\theta$. First we state a lemma that follows directly from \cite[Lemma 3.2]{GTV97}.

\begin{lemma}\label{l:32} 
Let $b_0, \gamma, a, a_1, a_2$ satisfy
\begin{align}
b_0&>\frac 12,\label{c0}\\
 0&\leq \gamma\leq 1,\label{c1}\\
 a, a_1, a_2&\geq 0,\label{c2}\\
 (1-\gamma)\max(a, a_1, a_2)&\leq b_0\leq (1-\gamma)(a+a_1+a_2),\label{c3}\\
 (1-\gamma) a&<b_0,\label{c4}\\
 m\geq \frac 52-(1-\gamma)(a+a_1+a_2)/b_0&\geq 0,\label{c5}
 \end{align}
with strict inequality\footnote{We will not need to consider this case.} in (\ref{c5}L) if equality holds in (\ref{c3}R) or if $a_1=0$.  And if
\begin{align}
\frac 12> \gamma a, \gamma a_1, \gamma a_2\label{c6},
\end{align}  
and $v, v_1, v_2 \in L^2(\R^3)$ are such that $\mathcal F^{-1} (\tfrac{ \hat v}{\ip {\tau+\abs{\xi}}^{\gamma a}}), \mathcal F^{-1} (\tfrac { \hat v_i}{\langle \tau_i+\abs{\xi_i}^2 \rangle^{\gamma a_i}}), i=1, 2, $ have support in $\abs{t}\leq CT$, then 
\begin{align}
\int\frac{\abs{\hat v \hat v_1\hat v_2}}{\ip{\tau+\abs{\xi}}^a\ip{\tau_1+\abs{\xi_1}^2}^{a_1}\ip{\tau_2+\abs{\xi_2}^2}^{a_2}\ip{\xi}^m }\lesssim T^\theta\norm{v}_2\norm{v_1}_2\norm{v_2}_2\label{m1a},\\
\int\frac{\abs{\hat v \hat v_1\hat v_2}}{\ip{\tau+\abs{\xi}}^a\ip{\tau_1+\abs{\xi_1}^2}^{a_1}\ip{\tau_2+\abs{\xi_2}^2}^{a_2}\ip{\xi_2}^m }\lesssim T^\theta\norm{v}_2\norm{v_1}_2\norm{v_2}_2,\label{m1b}
\end{align}
where
\be\label{t}
\theta=\gamma(a+a_1+a_2).
\ee
\end{lemma}

 \begin{remark} Lemma \ref{l:32} is a three dimensional version of  \cite[Lemma 3.2]{GTV97}, and the only difference is that $\frac 12$ upper bound in \eqref{c6} does not appear in \cite{GTV97}.  We included it here to maximize the value of $\theta$ appearing in \eqref{t}.  It also results in a simpler formula in \eqref{t}.  For the general formula for $\theta$ see \cite{GTV97}[(3.24) and (3.14)].   
 \end{remark}

We start by estimating $N_1$ and turn our attention to $N_2$ later.

 \begin{lemma}\label{l34} Let $0\leq \ell \leq 1$ and let $\eps_{0}, \eps,\bar{\eps}>0$ be sufficiently small.  Suppose the functions $v, v_1, v_2$ satisfy the conditions of Lemma \ref{l:32} and $b=\frac 12+{\bar{\epsilon}}$ and $c=1-\epsilon-b$, 
 then the estimate \eqref{m1} holds with
 \be\label{t1}
\theta=b+1-\epsilon-(\frac 52-\ell)(\frac 12+\epsilon_0).
 \ee
 \end{lemma}
 The proof is an application of Lemma \ref{l:32} and follows  \cite[Lemma 3.4]{GTV97}, but again, we attempt to maximize $\theta$. 
 \begin{proof}  Consider  two regions
 \begin{align*}
 \mbox{Region 1:}\quad \{\abs{\xi_1}\leq 2\abs{\xi_2} \},\\
 \mbox{Region 2:}\quad \{  \abs{\xi_1}>2\abs{\xi_2}  \}.
 \end{align*}
 In Region 1, \eqref{m1} reduces to
 \[
 \int\frac{\abs{\hat v \hat v_1\hat v_2}}{\ip{\tau+\abs{\xi}}^b\ip{\tau_1+\abs{\xi_1}^2}^{c}\ip{\tau_2+\abs{\xi_2}^2}^{b}\ip{\xi}^\ell }\lesssim T^\theta\norm{v}_2\norm{v_1}_2\norm{v_2}_2,
 \]
 which is exactly \eqref{m1a} with 
 \[
 (a, a_1, a_2, m)=(b, c, b, \ell).
 \]
 Therefore the estimate follows if we can find $0\leq\gamma\leq1$ and $b_0$ such that the conditions \eqref{c0}-\eqref{c6} are satisfied.  Let $b_0=\frac 12+\eps_0$.  One can check that if $\gamma=1-(\frac 52-\ell)\frac{b_0}{2b+c}$, then the conditions are satisfied, and we obtain 
  \[
 \theta=\gamma (2b+c)=b+1-\eps-(\frac 52-\ell)b_0,
 \]
 as needed.  
 
 Now we consider Region 2. Here we use \eqref{bound} to bound $N_1$ as follows
 \[
 N_1\lesssim I+I_1+I_2, 
 \]
 where 
 \begin{align*}
 I&=\int\frac{\abs{\hat v \hat v_1\hat v_2}}{\ip{\tau+\abs{\xi}}^{b-\frac{k-\ell}{2}}\ip{\tau_1+\abs{\xi_1}^2}^{c}\ip{\tau_2+\abs{\xi_2}^2}^{b}\ip{\xi_2}^k } ,\\
 I_1&=\int\frac{\abs{\hat v \hat v_1\hat v_2}}{\ip{\tau+\abs{\xi}}^{b}\ip{\tau_1+\abs{\xi_1}^2}^{c-\frac{k-\ell}{2}}\ip{\tau_2+\abs{\xi_2}^2}^{b}\ip{\xi_2}^k } ,\\
 I_2&=\int\frac{\abs{\hat v \hat v_1\hat v_2}}{\ip{\tau+\abs{\xi}}^{b}\ip{\tau_1+\abs{\xi_1}^2}^{c}\ip{\tau_2+\abs{\xi_2}^2}^{b-\frac{k-\ell}{2}}\ip{\xi_2}^k }.
 \end{align*}
 We apply estimate \eqref{m1b} of  Lemma \ref{l:32} three times. 
 Each time, we use the same value of $\gamma'$,  $0\leq \gamma'\leq 1$, which is chosen 
 so that  the resulting $\theta$ is the same in both regions.  This means
 \[
 \gamma'=\gamma\frac{2b+c}{2b+c-\tfrac{k-\ell}{2}}=\gamma\frac{2b+c}{2b+c-\tfrac 14}.
 \]
We let
\begin{align*}
(a, a_1, a_2, m)&=(b-\frac 14, c, b, k)\\
(a, a_1, a_2, m)&=(b, c-\frac 14, b, k)\\
(a, a_1, a_2, m)&=(b, c, b-\frac 14, k),
\end{align*}
and again one readily checks the conditions hold.  In particular, in each of the above three cases, when we verify \eqref{c3}R, it reduces to requiring
\[
b_0\leq (1-\gamma)(c+2b)-\frac 14=(\frac 52-\ell)b_{0}-\frac 14,
\]
which holds if and only if $\ell\leq \frac 32-\frac {1}{4b_{0}}$.
 \end{proof}

We turn our attention to treating $N_2$.

 \begin{lemma}\label{l35}Let $0\leq \ell \leq 1$, and let $\eps_{0}, \bar{\eps}, \eps_{1}>0$ be sufficiently small.  Suppose the functions $v, v_1, v_2$ satisfy the conditions of Lemma \ref{l:32} and $b=\frac 12+\bar{\epsilon}$ and $c=1-\epsilon-b$, then the estimate \eqref{m2} holds with
\be\label{t2}
\theta=b+1-\eps-(\frac 52-\ell)b_0.
 \ee
 \end{lemma}
 Note, $\theta$ here is the same as $\theta$ in Lemma \ref{l34}. The general  idea of the proof is the same as in Lemma \ref{l34} (also compare with  \cite[Lemma 3.5]{GTV97}).  We include the details for completeness.

 \begin{proof}  The proof is done in three regions, but due to symmetry of the estimate it is enough to consider only Region 1 and Region 2:
  \begin{align*}
 \mbox{Region 1:}&\quad \left\{\frac{\abs{\xi_2}}{2}\leq \abs{\xi_1}\leq 2\abs{\xi_2} \right\},\\
 \mbox{Region 2:}&\quad \left\{  \abs{\xi_1}>2\abs{\xi_2}  \right\}.\\
 \mbox{Region 3:}&\quad \left\{  \abs{\xi_1}\leq\frac 12\abs{\xi_2}  \right\}.
 \end{align*}
 In Region 1, \eqref{m2} reduces to
 \[
 \int\frac{\abs{\hat v \hat v_1\hat v_2}}{\ip{\tau+\abs{\xi}}^c\ip{\tau_1+\abs{\xi_1}^2}^{b}\ip{\tau_2+\abs{\xi_2}^2}^{b}\ip{\xi}^{2k-(\ell+1)} }\lesssim T^\theta\norm{v}_2\norm{v_1}_2\norm{v_2}_2,
 \]
 which is exactly \eqref{m1a} with 
 \[
 (a, a_1, a_2, m)=(c, b, b, \ell),
 \]
since $k=l+\frac 12$.  Here, if we also let $\gamma=1-(\frac 52 -\ell)\frac{b_0}{c+2b}$, then the estimate follows with $\theta$ given by \eqref{t2}.

 In Region 2 we again use \eqref{bound} to obtain
 \[
 N_2\lesssim I+I_1+I_2, 
 \]
 where 
 \begin{align*}
 I&=\int\frac{\abs{\hat v \hat v_1\hat v_2}}{\ip{\tau+\abs{\xi}}^{c-\frac 14}\ip{\tau_1+\abs{\xi_1}^2}^{b}\ip{\tau_2+\abs{\xi_2}^2}^{b}\ip{\xi_2}^k } ,\\
 I_1&=\int\frac{\abs{\hat v \hat v_1\hat v_2}}{\ip{\tau+\abs{\xi}}^{c}\ip{\tau_1+\abs{\xi_1}^2}^{b-\frac 14}\ip{\tau_2+\abs{\xi_2}^2}^{b}\ip{\xi_2}^k } ,\\
 I_2&=\int\frac{\abs{\hat v \hat v_1\hat v_2}}{\ip{\tau+\abs{\xi}}^{c}\ip{\tau_1+\abs{\xi_1}^2}^{b}\ip{\tau_2+\abs{\xi_2}^2}^{b-\frac 14 }\ip{\xi_2}^k }.
 \end{align*}
 As before, we apply Lemma \ref{l:32} and use estimate \eqref{m1b} three times with the same $0\leq \gamma'\leq 1$, which we choose so the resulting $\theta$ is the same in both regions.   Hence
  \[
 \gamma'=\gamma\frac{c+2b}{c+2b-\tfrac 14},
 \]
and
\begin{align*}
(a, a_1, a_2, m)&=(c-\frac 14, b, b, k)\\
(a, a_1, a_2, m)&=(c, b-\frac 14, b, k)\\
(a, a_1, a_2, m)&=(c, b, b-\frac 14, k).
\end{align*}
 
 \end{proof}
\section{Lower bound for the rate of blow-up of singular solutions}

We first summarize a priori estimates derived in previous sections.  Since \eqref{m1}-\eqref{m2} imply \eqref{m10}-\eqref{m20} combining with  \eqref{dh1}-\eqref{dh2} we have

\begin{align*}
\norm{\psi}_{X_S^{k,b}}&\lesssim \norm{\psi_0}_{H^{k}}+T^{\eps+\theta}\norm{\varphi_{2T}\mathcal Re~w}_{X_{W}^{\ell,b}}\norm{\varphi_{2T}\psi}_{X_S^{k,b}},\\
\norm{w}_{X_{W}^{\ell,b}}&\lesssim \norm{w_{0}}_{H^{\ell}}+T^{\eps+\theta}\norm{\varphi_{2T} \psi}_{X_S^{k,b}}^{2} +T^{1-b}\norm{\varphi_{2T}\mathcal Re~w}_{X_{W}^{\ell,0}}.
\end{align*}

Next applying \eqref{l25} with $q=2$ we have in particular
\[
T^{1-b}\norm{\varphi_{2T}\mathcal Re~w}_{X_{W}^{\ell,0}}\leq T^{1-b}\norm{\varphi_{2T}\mathcal Re~w}_{X_{W}^{\ell,b}}\leq T^{3/2-2b}\norm{\mathcal Re~w}_{X_{W}^{\ell,b}},
\]
and hence
\begin{align}
& \|\psi\|_{X_S^{k,b}}   \lesssim  \| \psi_0 \|_{H^k}
+ T^{\theta_{\ell}}\|\mathcal Re~w\|_{X_{W}^{\ell,b}} 
 \|\psi\|_{X_S^{k,b}}
\\
& \|w\|_{X_{W}^{\ell,b}}   \lesssim  \|w_{0}  \|_{H^\ell}
+ T^{\theta_{\ell}}  \|\psi\|^2_{X_S^{k,b}}+T^{3/2-2b}\|w\|_{X_{W}^{\ell,b}},
\end{align}
where
\[
 \theta_{l}=\eps + c+2b-(\frac 52-\ell)(\frac 12+\epsilon_0)-2b+1.
 \]
With the choice of $c=1-b-\eps=\frac 12-\bar{\eps}-\eps$ this reduces to
\be\label{theta}
\theta_{l}=\frac 14+\frac{\ell}{2}-\frac 52\epsilon_0+\ell\epsilon_0-\bar{\epsilon}=(\frac 14+\frac{\ell}{2})^{-} .
\ee

The two inequalities are rewritten as 
\begin{align*}
& \|\psi\|_{X_S^{k,b}}  + \|w\|_{X_{W}^{\ell,b}} \lesssim 
 \| \psi_0 \|_{H^k} + \|w_{0}  \|_{H^\ell}\cr
&\qquad +T^{\theta_{\ell}} \big( \|\psi\|_{X_S^{k,b}}+
\|w\|_{X_{W}^{\ell,b}} \big)^2+T^{3/2-2b}\|w\|_{X_{W}^{\ell,b}},
\end{align*}
or by using $$\|w\|_{X_{W}^{\ell,b}}\sim \|n\|_{X_{W}^{\ell,b}}+\|n_{t}\|_{X_{W}^{\ell-1,b}}\quad\mbox{and}\quad\|w_{0}  \|_{H^\ell}\sim \|n_{0}  \|_{H^\ell}+  \|n_{1}  \|_{H^{\ell-1}}$$  
equivalently we have
\begin{align}
\|\psi\|_{X_S^{k,b}}  + \|n\|_{X_{W}^{\ell,b}}+\|n_{t}\|_{X_{W}^{\ell-1,b}}& \le 
 C\big(\| \psi_0 \|_{H^k} +\|n_{0}  \|_{H^\ell}+  \|n_{1}  \|_{H^{\ell-1}}\big)\cr
&\quad  +CT^{\theta_{\ell}} \big( \|\psi\|_{X_S^{k,b}}+
 \|n\|_{X_{W}^{\ell,b}}+\|n_{t}\|_{X_{W}^{\ell-1,b}} \big)^2\nonumber\\
&\quad+CT^{3/2-2b}\big (\|n\|_{X_{W}^{\ell,b}}+\|n_{t}\|_{X_{W}^{\ell-1,b}}\big)\label{main}.
\end{align}

Further, we can assume $T$ is small enough so that $CT^{3/2-2b} <1/2$.  Rearranging \eqref{main} we obtain
\begin{align}
& \|\psi\|_{X_S^{k,b}}  + \|n\|_{X_{W}^{\ell,b}}+\|n_{t}\|_{X_{W}^{\ell-1,b}} \le 
2C( \| \psi_0 \|_{H^k} +\|n_{0}  \|_{H^\ell}+  \|n_{1}  \|_{H^{\ell-1}})\cr
&\qquad +2C T^{\theta_{\ell}} \big( \|\psi\|_{X_S^{k,b}}+
\|n\|_{X_{W}^{\ell,b}}+\|n_{t}\|_{X_{W}^{\ell-1,b}} \big)^2.
\end{align}
Then denoting 
\begin{align*}
{\mathcal X}(M,T) &=\left. \{ (\psi, n) :   (\psi, n, n_{t})\right|_{t =0}  =(\psi_{0}, n_{0}, n_{1}),\\
&\qquad\qquad~ \|\psi\|_{X^{k,b_1}_{S}}+
\|n\|_{X^{\ell,b}_{W}}+ \|n_{t}\|_{X^{\ell-1,b}_{W}} \le M \},
\end{align*}
the local wellposedness theory is obtained by a contraction argument in 
${\mathcal X}$ if $ 2CT^{\theta_{\ell}} M<1$ and $CT^{3/2-2b} <1/2$.

We adapt arguments developed Cazenave and Weissler
\cite{W81}, \cite{CW90} to prove a lower bound on the rate of blow-up. Denote by $t^*$ the supremum of all $T>0$ for which there exists a solution $(\psi, n)$ of the Zakharov system satisfying
$$ 
\|\chi_{[0,T]} \psi\|_{X_S^{k,b}}  + \|\chi_{[0,T]} n\|_{X_{W}^{\ell,b}}+\|\chi_{[0,T]} n_{t}\|_{X_{W}^{\ell-1,b}} < \infty.
$$
The local well-posedness theory shows $t^* >0$ and for all $t \in [0,t^*)$
$$
\| \psi (t) \|_{H^k} +\|n (t)  \|_{H^\ell}+  \|n_{t} (t)  \|_{H^{\ell-1}}< \infty.
$$
By maximality of $t^*$, it is impossible that
$$
\| \psi (t) \|_{{L^\infty_{[0,t^*]}} H^k} +\|n (t)  \|_{{L^\infty_{[0,t^*]}} H^\ell}+  \|n_{t} (t)  \|_{{L^\infty_{[0,t^*]}} H^{\ell-1}} < \infty.
$$
Otherwise, the initial value problem at time $t^*$ with Cauchy data $\psi (t^*), n(t^*), n_t (t^*)$ would be well-defined and the local theory would provide an extension of the solution beyond $t^*$. Therefore, if $t^* < \infty$, blow-up occurs:
\begin{equation}
\| \psi(t) \|_{H^k} + \|n(t)  \|_{H^\ell}+ \|n_t(t)  \|_{H^{\ell-1}}
\to \infty , \quad  {\rm as } \quad  t \to t^*. 
\end{equation}

Consider $\psi(t), n(t), n_t(t) $ posed at some time $t \in [0, t^*)$.  If for some $M$ 
\begin{equation} 
C(\| \psi(t) \|_{H^k} + \|n(t)\|_{H^\ell}
+ \|n_t(t)\|_{H^{\ell-1}})
+C (T-t)^{\theta_{\ell}} M^2 \le M
\end{equation}
then $T<t^*$. Therefore, $\forall M>0$,  
\begin{equation} 
C( \| \psi(t) \|_{H^k} + \|n(t)  \|_{H^\ell}+ \|n_t(t)  \|_{H^{\ell-1}})
+C (t^*-t)^{\theta_{\ell}} M^2 > M.
\end{equation}
Choosing $M=2C(\| \psi(t)\|_{H^k}+\|n(t)\|_{H^\ell}+ \|n_t(t)  \|_{H^{\ell-1}})$
we have  
\begin{equation}
M>C(t^*-t)^{-\theta_{\ell}},
\end{equation}
equivalently,
\begin{equation}
\| \psi(t)\|_{H^{\ell+1/2}}+\|n(t)\|_{H^\ell}+ \|n_t(t)\|_{H^{\ell-1}}>C(t^*-t)^{-\frac{1}{4}
(1+ 2\ell)^-}.
\end{equation}
which concludes the proof of  of  Theorem \ref{bup}.

\end{document}